\documentclass{amsart}
\usepackage{amsfonts}
\usepackage{amssymb}
\newtheorem{prop}{Proposition}

\newtheorem{cor}{Corollary}
\newtheorem{lem}{Lemma}
\begin{document}


\centerline{Appeared in: Differential Geometry and Applications, 7 (1997), 123-130}

\vspace{1cm}

\title[Focus-focus singularities]{A note on focus-focus singularities}

\author{Nguyen Tien Zung}
\address{SISSA - ISAS, Via Beirut 2-4, Trieste 34013, Italy (old address)}
\email{tienzung@math.univ-montp2.fr}
\date{27 May 1994, revised 29 March 1995}
\keywords{integrable Hamiltonian system, singularity, monodromy.}
\subjclass{58F07,70H05}

\begin{abstract} {
We give a topological and geometrical description of focus-focus singularities of
integrable Hamiltonian systems. In particular, we explain why the monodromy around
these singularities is non-trivial, a result obtained before by J.J. Duistermaat and
others for some concrete systems. }
\end{abstract}
\maketitle

\section{Introduction}

Many integrable Hamiltonian systems in classical mechanics - from as simple as the
spherical pendulum - contain focus-focus singularities (see Section
\ref{section:example}). Thus the study of these singularities is important in order
to understand the topology of integrable systems. We address this problem in the
present note. It turns out that the topological structure of focus-focus
singularities is quite simple, though very different from elliptic and hyperbolic
cases. The affine structure of the orbit space near focus-focus singularities is
also very simple (cf. Proposition \ref{prop:orbitspace}). As a corollary, we obtain
that the monodromy around these singularities is non-trivial. The notion of
monodromy was first given by Duistermaat \cite{Duistermaat}, and its non-triviality
was observed by Duistermaat, Cushman, Kn\"orrer, Bates, etc., for various systems,
all of which turn out to be connected with focus-focus singularities (see Sections
\ref{section:orbitspace},\ref{section:example}).

For simplicity of the exposition we will
consider only systems with two degrees of freedom.
The results remain unchanged for focus-focus codimension 2 singularities
of integrable Hamiltonian systems with more degrees of freedom.

\section{Local structure}

Throughout this work, by an integrable system we will mean a Poisson ${\mathbb R}^2$
action on a real smooth symplectic 4-manifold $(M^4, \omega)$, given by a moment map
${\bf F} = (F_1, F_2) : M^4 \to {\mathbb R}^2$. We will also assume that the level
sets of ${\bf F}$ are compact, and hence they are disjoint unions of Liouville tori
wherever non-singular.

Suppose that $x_0 \in M^4$ is a fixed point of the above Poisson action:
$dF_1 (x_0) = dF_2 (x_0) = 0$. Let $H_i$ be the quadratic part of $F_i$ at
$x_0$ $(i = 1,2)$. Since $F_1, F_2$ are Poisson commuting, so are $H_1$
and $H_2$: $\{ H_1, H_2 \} = 0$.

We will assume that $x_0$ is a non-degenerate singular point, i.e. $H_1$ and
$H_2$ form a Cartan subalgebra of the algebra of quadratic forms under the
natural Poisson bracket. Then by a classical theorem of
Williamson (see e.g. \cite{Williamson,Arnold,LU}), after a linear change of
basis of the
Poisson action (i.e. a linear change of the moment map ${\bf F}' = A \circ
{\bf F}$, $A$ being a constant invertible matrix), one of the following four
alternative cases happens:

$H_1 = x_1^2 + y_1^2, H_2 = x_2^2 + y_2^2,$  (elliptic case)

$H_1 = x_1^2 + y_1^2, H_2 = x_2 y_2,$ (elliptic-hyperbolic case)

$H_1 = x_1y_1, H_2 = x_2 y_2,$  (hyperbolic case)

$H_1 = x_1 y_1 + x_2 y_2,  H_2 = x_1 y_2 - x_2 y_1.$ (focus-focus case)

Here $(x_1, y_1, x_2, y_2)$ is a system of symplectic coordinates in the
tangent space at $x_0$: $\omega_{x_0} = dx_1 \wedge dy_1 + dx_2 \wedge
dy_2$.

In this paper we are interested in the fourth, focus-focus case. Notice that if one
considers systems with complex coefficients, then the above four cases become the
same. It is equivalent to say that the complex symplectic algebra $sp(2n, {\mathbb
C})$ has only one conjugacy class of Cartan subalgebras.

The local analysis of focus-focus singular points was done by Lerman and Umanskii
\cite{LU} and Eliasson \cite{Eliasson}.
It turns out that near a focus-focus singular point $x_0$,
there are two local Lagrangian invariant submanifolds which intersect
transversally at $x_0$ (cf. \cite{LU}), and the other nearby local invariant
Lagrangian submanifolds are annuli. Moreover, the local singular Lagrangian
foliation given by the moment map ${\bf F}$ near $x_0$ is {\it
symplectically}
equivalent to the one given by the linearized moment map ${\bf H} = (H_1,
H_2)$ (cf. \cite{Eliasson}).
In other words, there is a local symplectic system of
coordinates $(x_1, y_1, x_2, y_2)$ around $x_0$, for which ${\bf F}$ can be
expressed as a function of two variables $f_1, f_2$, where $f_1 = x_1y_1 +
x_2 y_2, f_2 = x_1y_2 - x_2y_1$. Such a local system of symplectic
coordinates will be called a {\it canonical system of coordinates} near a
focus-focus point $x_0$.

\begin{prop}
Let $x_0$ be a focus-focus singular point. Then there is a natural
(local) circle Hamiltonian action near $x_0$, which preserves the moment map.
\end{prop}

{\it Proof}. Let $(x_1, y_1, x_2, y_2)$ be a canonical system of coordinates. Then
the function  $f_2 = x_1y_2 - x_2y_1$ will be the required Hamiltonian. Remark that
this Hamiltonian ${\bf S^1}$ action has exactly one fixed point, namely $x_0$.
$\square$

In fact, one can give a different proof of the above proposition, which relies only
on the existence of two transversal invariant Lagrangian submanifolds indicated
above, as follows: by perturbation theory, one sees that there is a linear
combination of the Hamiltonian vector fields $X_{F_1}, X_{F_2}$, which gives rise to
a periodic flow on one of these Lagrangian submanifolds. Then one can extend this
flow in the most natural way to obtain the required Hamiltonian action in a
neighborhood of $x_0$. Notice also the uniqueness of this Hamiltonian ${\bf S}^1$
action, up to the direction.

The following proposition is not needed for the rest of this note but will be used
to obtain some geometric invariants of torus singular Lagrangian foliations with
focus-focus singularities.

\begin{prop}
\label{prop:canonical}
Let $(X_1, Y_1, X_2, Y_2)$ and $(x_1, y_1, x_2, y_2)$ be two different
canonical systems of coordinates at
$x_0$. Then we have $X_1Y_2 - X_2 Y_1 = x_1 y_2 - x_2 y_1$ or
$X_1Y_2 - X_2 Y_1 = - (x_1 y_2 - x_2 y_1)$, and the
difference between $X_1Y_1 + X_2Y_2$ and $\pm (x_1y_1 + x_2y_2)$  is a
flat function.
\end{prop}

{\it Proof}. The proof follows directly from the results of Vey \cite{Vey} (who,
however, seems to omit the focus-focus case). In the analytic case, one can
complexify the system so that to return the problem to the case of elliptic
singularity. Then one sees that $f_1 = x_1 y_1 + x_2 y_2$ and $f_2 = x_1y_2 -
x_2y_1$ are some linear combinations (namely, the sum and the difference, up to some
scalars) of two action functions (i.e. action components of a local system of
action-angle coordinates near an elliptic singularity). Hence $f_1$ and $f_2$ do not
depend on a particular choice of canonical coordinates, up to a sign (cf.
\cite{Vey}). In the smooth case, Taylor expansions will give the same result, up to
flat functions. Note that even in the smooth case we have $X_1Y_2 - X_2 Y_1 = x_1
y_2 - x_2 y_1$ or $X_1Y_2 - X_2 Y_1 = - (x_1 y_2 - x_2 y_1)$ since they are the
Hamiltonian of a unique natural ${\bf S}^1$ action discussed in the previous
proposition. $\square$

\section{Stable case: topology and monodromy}
\label{section:orbitspace}

Denote by $N(x_0)$ the connected component of the preimage of the moment map
${\bf F}$ which contains $x_0$. We will always assume that all singular
points in $N(x_0)$ of the Poisson action are non-degenerate.
(See e.g. \cite{Eliasson,LU,Zung} for the
definition of nondegeneracy). Then $N(x_0)$ is a
non-degenerate singular leaf in the singular Lagrangian foliation by
Liouville tori in a most natural sense (see \cite{Zung}
for more details). From the
results of Lerman and Umanskii \cite{LU} it follows that singular points in
$N(x_0)$ either lie in one-dimensional
closed singular hyperbolic orbits or are
focus-focus fixed points.

By convention, we will say that a focus-focus singular leaf $N(x_0)$ is {\it
topologically stable} if it does not contain singular hyperbolic orbits, i.e. if all
singular points in it are focus-focus fixed points.

Suppose now that $N(x_0)$ is topologically stable and contains exactly $n$
focus-focus fixed points $x_0,\ldots, x_{n-1}$. Then because of the Poisson
${\mathbb R}^2$ action, $N(x_0) \setminus \{ x_0, \ldots, x_{n-1} \} $ must be a
non-empty disjoint union of annuli. It follows that $N(x_0)$ consists of a chain of
$n$ Lagrangian spheres, each of which intersects transversally with two other. (This
simple but important fact was observed by Bolsinov, and also by Lerman and Umanskii
themselves). In particular, the fundamental group of a tubular neighborhood of it is
${\mathbb Z}$. When $n=1$, $N(x_0)$ is just a sphere with one point of
self-intersection. It is well-known that the orbit space of the singular Lagrangian
foliation (by Liouville tori) has a unique natural integral affine structure outside
the singularities (see e.g. \cite{Duistermaat}). We have the following:

\begin{prop}
\label{prop:orbitspace} Let $N(x_0)$ be a topologically stable non-degenerate
focus-focus leaf with $n$ fixed points as above. Then in a neighborhood of the image
of this leaf in the orbit space of the singular Lagrangian foliation, the affine
structure can be obtained from the standard flat structure in ${\mathbb R}^2 = \{(x,
y)\}$ near the origin $O$ by cutting out the angle $\angle\{(0,1),(-n,1)\}$ and
gluing the edges of the rest together by the integral linear transformation $(x, y)
\mapsto (x + ny, y)$.
\end{prop}

Before proving the above proposition let us now construct an algebraic model for
this singularity.

Near the origin $O$ in the local standard symplectic space ${\mathbb R}^4, \omega  =
dp_1 \wedge dq_1 + dp_2 \wedge dq_2 $ we have two generating functions for a Poisson
${\mathbb R}^2$-action with the singularity of the type focus-focus:
$$ f_1 = p_1 q_1 + p_2 q_2, f_2 = p_1 q_2 - p_2 q_1$$
Set $z_1 = p_1 - i p_2 , z_2 = q_1 + iq_2 $ . Here
 we define a complex structure:
$$ J: (p_1,p_2,q_1,q_2) \to (p_2,-p_1,-q_2,q_1) $$
Then $f_1$ and $f_2$ are the real and imaginary part of function $z_1 z_2$. In
particular, the level sets of $(f_1,f_2)$ are the level sets  of $z_1z_2$ in
${\mathbb C}^2$.

Consider first the simplest case, when $n=1$. We can construct the model as follows:
Take the conformal map $\Phi: (z_1,z_2) \mapsto (z_2^{-1}, z_1z_2^2)$. Note that
$\omega = Re dz_1 \wedge dz_2$, and $\Phi$ is a complexification of a real
area-preserving map. Hence $\Phi$, where it is well-defined, is a symplectic
mapping. Consider a small neighborhood $ D \times {\mathbb C}P^1$ of the sphere $0
\times {\mathbb C}P^1$ ($z_2$ lies in ${\mathbb C}P^1$). Gluing the points near
(0,0) to the points near $(\infty, 0)$ by the map $\Phi$ we obtain a complex space
${\mathcal U}$ which has a natural symplectic form $\omega$ (because $\Phi$
preserves the symplectic form). Furthermore, it is
 clear that the analytic map $z_1z_2: {\mathcal U} \to {\mathbb C}$ is
well-defined on M and is the moment map for a desired ${\mathbb R}^2$ Poisson
action.

The case $n > 1$ is similar. Take $n$ samples $U_i (i=1,...,n)$ of $ D \times
{\mathbb C}P^1$ and define $n$ local maps $\phi_i : U_i \to U_{i+1} (U_{n+1} =
U_1)$, which in local coordinates have the same form $\Phi$ as above. Glue $U_i$
together by these maps. Again we obtain a symplectic manifold, and the function $z_1
z_2$ provides us a moment map for a Poisson ${\mathbb R}^2$ action on that manifold,
now with $n$ focus-focus points over a  singular point of the local bifurcation
diagram. More geometrically, what we do is just take an $n$-covering of the
one-point case. Under this covering every Lagrangian torus also pulls back to an
$n$-covering of itself.

Topologically, this construction is unique, i.e. it can be easily shown that any two
topologically stable focus-focus singularities $N_1(x_0)$ and $N_2(y_0)$ with the
same number of singular points will have diffeomorphic\footnote{erratum:
homeomorphic, not necessarily diffeomorphic} singular Lagrangian foliations (though
in general we don't have a foliation-preserving {\it symplectomorphism} between
them, cf. \cite{Grossi}). If we forget about the foliation, then by Moser's path
method, one can show that there is a symplectomorphism between some tubular
neighborhoods ${\mathcal U}_1(N_1)$ and ${\mathcal U}_2(N_2)$ of $N_1$ and $N_2$,
which sends $N_1$ to $N_2$ (cf. \cite{Moser,Weinstein,Brunella}).

From the construction it is easy to see the topological type of  an
``isoenergy''
3-manifold around the singularity (more precisely, the manifold $\{ | z_1z_2
| = \epsilon > 0 \} $). It is a locally flat fibration with torus fiber
over a
circle.  We will compute the holonomy mapping of this fiber bundle. Our
computations in fact do not depend on the above specific model, but only on
topological properties of all topologically stable focus-focus singularities.

Recall from the previous section that locally near $x_0$ there is a natural
Hamiltonian ${\bf S}^1$ action. This action can be extended to be a Hamiltonian
${\bf S}^1$ action in a tubular neighborhood ${\mathcal U}(N)$ of $N(x_0)$. Denote
by $g$ the corresponding Hamiltonian function, $g(x_0) = 0$. Because of the
invariance, $g$ can also be considered as a function on the orbit space. In our
model $g = f_2$. Fix a small circle $\{|z_1 z_2| = \epsilon \}= \{ f_1^2 +f_2^2 =
\epsilon^2 \}$ in the orbit space. Every point in this circle corresponds to one
Liouville torus. Fix one point $\{ f_1=\epsilon,f_2=0\}$. On the torus corresponding
to this point fix a basis of generators of the fundamental group, so that the second
generator is induced from the symplectic vector field $X_{f_2}$, and the orientation
on ${\bf T}^2$ given by these two generators coincides with that one given by
$X_{f_1}, X_{f_2}$. Denote these generators by $\gamma, \delta$ respectively.  When
the point $\{ f_1=\epsilon,f_2=0\}$ moves along the circle in the positive direction
(anti-clockwise), $\gamma$ and $\delta$   also move homotopically, and in the end
come back to some new cycles $\gamma_{new} , \delta_{new}$ on the old torus.

\begin{lem}
\label{lem:focus}
With the above notations we have:\\
1) The curves $\{g = {\rm const}\}$ are straight
lines in the affine structured
orbit space. \\
2) $\begin{pmatrix} \gamma_{new} \cr \delta_{new} \end{pmatrix} = \begin{pmatrix} 1
& n \cr 0 & 1 \end{pmatrix} \begin{pmatrix} \gamma \cr \delta \end{pmatrix} $
\end{lem}

{\it Proof}. 1) follows from the fact that  the flow of $X_g$ is periodic with
constant period. We prove 2) for $n=1$. After that one can use the n-covering
argument to see it for any n. Since the assertion is topological, then it is enough
to prove it in our model. On the submanifold    $\{|z_1 z_2| = \epsilon \}$ in
${\mathcal U}$ denote by $\theta$ the cycle where $z_1 = {\rm const}$ and arg $z_2$
decreases, $\lambda$  the circle where $z_2 = {\rm const}$ and arg $z_1$ increases.
Then $\delta = \theta + \lambda$. Recall that when we go around by $\gamma$, the
coordinate system changes by the rule: $ (z_1^{new}, z_2^{new}) = (z_2^{-1},
z_1z_2^2)$. Fix a point ${\mathcal A}$ in one representative of $\gamma$. Moving
${\mathcal A}$ along $\lambda$ by some angle $e$ means increasing arg $ z_1 $ by
this angle. Making things go homotopically around $\gamma$, what we get is that
arg$z_1^{new}$ increases by $e$, $z_2^{new}$ remains constant. By the above rule, in
the old coordinates arg $z_1$ increases by 2$e$, and arg $z_2$ decreases by $e$.
That yields that after going around $\gamma$, ${\mathcal A}$ becomes to move on
$\lambda + \delta$ with the same angle. It follows that $\gamma_{new} = \gamma +
\delta$. $\square$

\vspace{0.5cm} {\it Proof of Proposition~\ref{prop:orbitspace}}. It follows directly
from Lemma~\ref{lem:focus}. $\square$ \vspace{0.5cm}

Let $V^m$ be an affine structured manifold. Then in the tangent bundle of
$V$ there is a unique natural flat connection, and fixing a point $x \in V$
there is a monodromy linear representation of $\pi_1(V)$ in $T_xV$, defined
as usual (cf. \cite{Duistermaat}).
From Proposition \ref{prop:orbitspace} we immediately get:

\begin{cor}
The local monodromy near every topologically stable focus-focus point in the orbit
space is nontrivial (and is generated by $\begin{pmatrix} 1 & n \cr 0 & 1
\end{pmatrix}$).
\end{cor}

\section{${\bf S}^1$ action and reduction}
\label{section:reduction}

Consider now the non-stable case as well, i.e. allow $N(x_0)$ to contain hyperbolic
singular orbits. As before, 2-dimensional orbits in $N(x_0)$ are annuli. It follows
that $N(x_0)$ is a union of immersed closed Lagrangian surfaces which intersect
transversally at hyperbolic orbits and focus-focus points. Again, it can be easily
seen that the local ${\bf S}^1$ action discussed before can be extended naturally to
a Hamiltonian ${\bf S}^1$ action in a saturated neighborhood ${\mathcal U}(N(x_0))$
of $N(x_0)$, which preserves the moment map. Near (possible) hyperbolic orbits in
$N(x_0)$, this ${\bf S}^1$ action coincides with another Hamiltonian ${\bf S}^1$
action, which is defined in a natural way in a tubular neighborhood of each closed
hyperbolic orbit of the Poisson action. Notice that the natural ${\bf S}^1$ action
defined near a hyperbolic closed orbit has isotropy group at most ${\mathbb Z}_2$
(the cyclic group of two elements) at this and some nearby hyperbolic orbits and is
free outside them (cf. \cite{Zung}). In other words, we have:

\begin{prop}
i) In a saturated neighborhood  ${\mathcal U}(N(x_0))$ there is a unique Hamiltonian
${\bf S}^1$ action, generated by a function $g$, $g(x_0) = 0$, which preserves the
moment map. In particular, it leaves $N(x_0)$ and hyperbolic
singular orbits invariant. \\
ii) This action is trivial at focus-focus points, may have isotropy group ${\mathbb
Z}_2$ at hyperbolic orbits, and is free elsewhere. $\square$
\end{prop}

Consider the moment map $g: {\mathcal U}(N) \to {\mathbb R}$ of the above ${\bf
S}^1$ action. At each small value $s$ denote by $P_s$ the symplectic 2-dimensional
space obtained by the Marsden-Weinstein reduction at $g = s$.

Consider $P_0$. It contains the image of $x_i$, denoted by $p_i$ (i = 0,...,n-1).
Let $(x_1, y_1, x_2, y_2)$ be a canonical system of coordinates at $x_0$. Then each
orbit of the ${\bf S}^1$ action, which lies in $\{ g = 0 \}$, intersects the
symplectic plane $x_1 = y_1 = 0$, and the intersection is a pair of points of the
type $\{ (x_2, y_2), (-x_2, -y_2) \}$. It follows that $P_0$ is an orbifold of order
2 at $p_i$. (Homeomorphically we can `desingularize' the points $p_i$, but not
symplectically). Let $q_1,...,q_k, k \geq 0$ denote the image of
normally-nonorientable hyperbolic orbits (i.e. orbits on which the ${\bf S}^1$
action is not free) of $N(x_0)$ in $P_0$. Then $P_0$ is also an orbifold of order 2
at these points. Thus, $P_0$ is a topological surface, but symplectically it is a
quotient of a symplectic surface by a ${\mathbb Z}_2$ action.

Notice that, since $g$ can be viewed as a function on the orbit space of the
original ${\mathbb R}^2$ action, the restriction of the moment map on $\{ g=s\}$
will give rise to a circle foliation on $P_s$. On $P_0$ this foliation is singular,
with the singular leaf being the image of $N(x_0)$. In the topologically stable case
this singular leaf is just a circle which contains all of the points $p_i$. In the
non-stable case, $P_0$ with the singular foliation looks like a hyperbolic
codimension 1 singularity (cf. \cite{Zung}), only it contains some special
(focus-focus) points in the singular leaf. In the topologically stable case, $P_s$
($s \neq 0$) with the circle foliation on it is regular. In the non-stable case it
can be obtained from $P_0$ by smoothening the points $p_i$ and perturbing the
foliation a little bit. In particular, each $P_s, s \neq 0$ now represents a
codimension 1 singularity.

The above description of the Marsden-Weinstein reduction for the
distinguished ${\bf S}^1$ action gives a better understanding of the
topology of topologically stable and non-stable focus-focus points. On the other hand, by
a small perturbation of the Poisson action we can always split out
focus-focus singularities from hyperbolic codimension 1 singularities. In
other words, we have:

\begin{prop}
If $N(x_0)$ contains hyperbolic orbits, then there is an arbitrarily $C^{\infty}$
small perturbation ${\bf F}'$ of ${\bf F}$, such that ${\bf F}'$ is again the moment
map of some Poisson ${\mathbb R}^2$ action, which has $x_0$ as a focus-focus point,
and the singular leaf $N'(x_0)$ with respect to ${\bf F}'$ contains the same number
of focus-focus points as $N(x_0)$ but does not contain any hyperbolic orbit.
\end{prop}

{\it Proof}. By a local diffeomorphism of ${\mathbb R}^2$, we can assume that ${\bf
F}(x_0) = 0$ and $F_2 = g$, i.e. it generates the distinguished Hamiltonian ${\bf
S}^1$ action.

If $\gamma$ is a closed 1-dimensional hyperbolic orbit in $N(x_0)$, and the ${\bf
S}^1$ action on $\gamma$ is free, then it is easy to construct a system of
symplectic coordinates $(x_1, y_1, x_2, y_2)$, $y_2$ - mod $1$, near $\gamma$, such
that $\gamma = \{ x_1 = x_2 = y_2 = 0\}$, and $x_2 = g$. In this canonical system of
coordinates, $F_1$ is a function depending only on 3 variables $x_1, x_2, y_2$.
Moreover, we can make so that in a small tubular neighborhood of $\gamma$, $N(x_0)$
is given by $ N(x_0) \cap {\mathcal U}(\gamma) = \{F_1 = F_2 = 0\} \cap {\mathcal
U}(\gamma) = \{ x_1 =0, x_2y_2 =0\}$. Then we can slightly perturb $F_1$, as a
function of three variables $(x_1, x_2, y_2)$, so that it remains unchanged outside
${\mathcal U}(\gamma)$, and $\{F_1 = F_2 = 0\} \cap {\mathcal U}(\gamma)$ becomes
smooth.

In case the ${\bf S}^1$ action is free on all hyperbolic orbits of $N(x_0)$, we can
apply the above procedure to all these hyperbolic orbits to obtain the required
result. In case there are some orbits with isotropy group ${\mathbb Z}_2$, we can
use a double covering and make everything ${\mathbb Z}_2$ invariant to obtain the
same result. $\square$

The above proposition gives a justification for the word {\it stable}.
Unlike the case of complicated hyperbolic (codimension 1) singularities, we
don't fear that when we perturb the things using the above proposition some
finite symmetry breaks up, since focus-focus points and hyperbolic orbits
are clearly of different natures and have different codimensions. One can
also `split' focus-focus points, i.e. make them lie on different levels
of the moment map, by a similar ${\bf S}^1$-invariant perturbation. But
then some good finite symmetry may break up.

\section{Examples and remarks}
\label{section:example}

We have given the topological classification, and the affine structure of the orbit
space, for focus-focus singularities. The geometrical classification (i.e. up to
foliation-preserving symplectomorphisms) is discussed in work \cite{Grossi}, where
it is shown that there arises some formal Taylor series in the set of invariants,
like in \cite{DMT}.

In \cite{Zung} we have shown that (topological) 2-domains of orbit spaces of
integrable Hamiltonian systems\footnote{erratum: the system must not contain some
types of degenerate singularities} with two degrees of freedom can have fundamental
group at most ${\mathbb Z}^2$, and 2-domains with non-trivial fundamental group can
appear only in very special systems. Thus in general, at least for systems with two
degrees of freedom, non-trivial monodromy is most probably connected with
focus-focus singularities. (Recall that topological 2-domains of orbit spaces can
contain focus-focus points).

It seems that the study of action-angle variables plays an important role in
classical mechanics (see e.g. a survey by Marle \cite{Marle}). A particular
attention is given to the `phenomenon' of nontriviality of the monodromy. As we
speculated above, this phenomenon is almost for sure connected to the existence of
topologically stable focus-focus singularities (but see \cite{Zung}). We list here
some known examples:

1. {\it Spherical pendulum} (cf. \cite{Duistermaat}). The spherical pendulum has a
${\bf S}^1$ group of symmetries (rotations), hence it is an integrable system with
two degrees of freedom. It has 2 stationary points: the lowest and highest positions
with zero velocity. The lowest position is stable, and indeed it is an elliptic
singular point. The highest position is unstable dynamically, and one can see that
it is a focus-focus singular point, by just looking at the trajectories having this
position as the limit. It follows that we have a stable (in our topological sense)
focus-focus singularity with one fixed point.

2. {\it Lagrange top}. A detailed analysis of this classical spinning top is given
in \cite{CK}, together with the nontrivial monodromy. The existence of a focus-focus
singularity was also observed by many people (see e.g. \cite{Eliasson,Oshemkov}).

3. {\it Champagne bottle} (cf. \cite{Bates}). The Hamiltonian is $\frac{1}{2} (p_x^2
+ p_y^2) - (x^2 + y^2) + (x^2 + y^2)^2$, i.e. a special case of Garnier systems.
Bates showed that there is a focus-focus singularity, and the monodromy is generated
by the same matrix $\begin{pmatrix} 1 & 0 \cr 1 & 1 \end{pmatrix}$ as in the
previous examples (or $\begin{pmatrix} 1 & 1 \cr 0 & 1 \end{pmatrix}$ if one
permutes the basis).

4. {\it Clebsch's equation} (motion of a rigid body in a fluid). The bifurcation
diagram of this system was constructed by Pogosyan \cite{Pogosyan}, from where the
existence of a focus-focus singularity is clear.

5. {\it Euler's equation on $so(4)$}. The bifurcation diagram of some integrable
Euler's equations in $so(4)$ was constructed by Oshemkov (see, e.g.,
\cite{Oshemkov}). These bifurcation diagrams also contain some isolated singular
points, i.e. focus-focus points! One can suspect that Euler's equations in many
other Lie (co)algebras will also possess focus-focus singularities.

After this note was written, I found two relevant papers \cite{LU2} and \cite{Zou}.
Lerman and Umanskii \cite{LU2} also studied the topology of extended neighborhoods
of focus-focus singularities, but their description is rather complicated. Zou
\cite{Zou} already proved Corollary 1, but only for the case $n=1$. His proof is
based on an interesting observation that in case $n=1$, the situation resembles the
simplest case of Picard-Lefschetz theory. (In fact, our model in Section 3 is
holomorphic so one can apply Picard-Lefschetz theory). In \cite{Zou} Zou also
mentioned a focus-focus singularity with $n=2$ (in a system studied by him and Larry
Bates), where his theorem does not apply.

\vspace{0.5cm} {\bf Acknowledgements}. This paper is based on a part of my thesis
\cite{Zung}. I would like to thank my supervisor Mich\`ele Audin for her enormous
help, and Pierre Molino for pointing out to me the importance of Duistermaat's work
\cite{Duistermaat}, and the reference \cite{Marle}. I would like also to thank the
referee for his critical remarks.

\bibliographystyle{amsplain}

\end{document}